\documentclass{amsart}

\usepackage{amsmath}
\usepackage{amsfonts}
\usepackage{amssymb}

\numberwithin{equation}{section}

\newtheorem{theorem}[subsection]{Theorem}
\newtheorem{proposition}[subsection]{Proposition}
\newtheorem{lemma}[subsection]{Lemma}

\newcommand{\wh}{\widehat}

\begin{document}

\title{A note on Fre{\u\i}man's theorem in vector spaces}

\author{Tom Sanders}
\address{Department of Pure Mathematics and Mathematical Statistics\\
University of Cambridge\\
Wilberforce Road\\
Cambridge CB3 0WA\\
England } \email{t.sanders@dpmms.cam.ac.uk}

\begin{abstract}
A famous result of Fre{\u\i}man describes the sets $A$, of integers,
for which $|A+A| \leq K|A|$. In this short note we address the
analagous question for subsets of vector spaces over
$\mathbb{F}_2$. Specifically we show that if $A$ is a subset of a
vector space over $\mathbb{F}_2$ with $|A+A| \leq K|A|$ then $A$
is contained in a coset of size at most $2^{O(K^{3/2}\log K)}|A|$,
which improves upon the previous best, due to Green and Ruzsa, of
$2^{O(K^2)}|A|$. A simple example shows that the size may need to
be at least $2^{\Omega(K)}|A|$.
\end{abstract}

\maketitle

\section{Introduction}

If $A$ and $B$ are subsets of an abelian group $G$ then we define
the \emph{sumset}, $A+B$, to be the set of all elements formed by
adding an element of $A$ to an element of $B$ i.e. $A+B=\{a+b:a
\in A \textrm{ and }b \in B\}$. There is a famous result of
Fre{\u\i}man \cite{GAF} which in some sense describes the sets $A
\subset \mathbb{Z}$ for which $|A+A|\leq K|A|$. This note concerns
what in modern parlance would be called the finite field analogue
of Fre{\u\i}man's result. Specifically it concerns the following
theorem.
\begin{theorem}\label{FfF}\emph{(Finite Field Fre{\u\i}man)}
Suppose that $G$ is a vector space over $\mathbb{F}_2$. Suppose
that $A \subset G$ is a finite set with $|A+A| \leq K|A|$. Then
$A$ is contained in a coset of size at most $f(K)|A|$.
\end{theorem}
While finite field models are an important tool for understanding
problems in general abelian groups this result has independent
significance in coding theory and has been pursued by a number of
authors. We do not attempt a comprehensive survey here, but
mention a few papers which are important from our standpoint.

The paper \cite{JMDFHAP} of Deshouillers, Hennecart and Plagne
provides an overview of the problem and records a quantitative
version of Theorem \ref{FfF} due to Ruzsa; it is a relatively
simple argument which shows that one may take $f(K) \leq
K2^{\lfloor K^3\rfloor -1}$. The bulk of their paper concerns
refined estimates for the case when $K$ is small; by contrast our
interest lies in the asymptotics.

In a recent paper of Green and Ruzsa \cite{BJGIZR1} the authors
improve Ruzsa's bound from \cite{JMDFHAP} when they show that one may
take $f(K) \leq K^22^{\lfloor 2K^2-2\rfloor}$. In this note we
refine this further by proving that one may take $f(K) \leq
2^{O(K^{3/2}\log K)}$. Formally, then, we shall prove the
following theorem.
\begin{theorem}\label{maintheorem}
Suppose that $G$ is a vector space over $\mathbb{F}_2$. Suppose
that $A \subset G$ is a finite set with $|A+A| \leq K|A|$. Then
$A$ is contained in a coset of size at most $2^{O(K^{3/2}\log
K)}|A|$.
\end{theorem}
For comparison we record the following well known example. Let $H$
be a finite subgroup of $G$ and $g_1+H,...,g_{K-1}+H$ be $K-1$
linearly independent cosets of $H$ in the quotient space $G/H$.
Let $A$ be the union of $H$ and the representatives
$g_1,...,g_{K-1}$. Then $|A|=|H| + K-1 \sim |H|$ and
\begin{equation*}
|A+A|=K(|H|+(K-1)/2) \sim K|H| \lesssim K|A|.
\end{equation*}
However $A$ contains a linearly independent set of size $\dim H +
K-1$ and so $A$ is not contained in a coset of dimension less than
$\dim H + K-2$ hence if $H'$ is a coset containing $A$ then
\begin{equation*}
|H'| \geq 2^{\dim H + K-2} = 2^{K-2}|H| \gtrsim 2^{K-2}|A|.
\end{equation*}
This example perhaps suggests that one could take $f(K) \leq
2^{O(K)}$ in Theorem \ref{FfF}. In fact there are other more
compelling reasons to believe this, however it does not seem to
reflect the underlying situation; in \cite{BJGPFR} Green addresses
this concern by introducing (a special case of) the Polynomial
Fre{\u\i}man-Ruzsa conjecture (attributed to Marton in \cite{IZRArb})
which, if true, seems to have some very important applications.
For a detailed discussion of this see either the paper of Green or
Chapter 5 of the book \cite{TCTVHV} of Tao and Vu.

\section{Proof of Theorem \ref{maintheorem}}

In \cite{BJGIZR}, Green and Ruzsa extended Ruzsa's proof of
Fre{\u\i}man's theorem from \cite{IZRF} to arbitrary abelian groups; for
an exposition see \cite{BJGEdin}. Our proof is a refinement of
theirs.

Their method becomes significantly simpler in the vector space
setting, and would immediately give us the following weak version
of the main theorem.
\begin{theorem}
Suppose that $G$ is a vector space over $\mathbb{F}_2$. Suppose
that $A \subset G$ is a finite set with $|A+A| \leq K|A|$. Then
$A$ is contained in a coset of size at most $2^{O(K^2\log K)}|A|$.
\end{theorem}
The proof involves three main step.
\begin{itemize}
\item \emph{(Finding a good model)} First we use the fact that
$|A+A| \leq K|A|$ to show that $A$ can be embedded as a
\emph{dense} subset of $\mathbb{F}_2^n$ in a way which preserves
much of its additive structure. \item \emph{(Bogolubov's
argument)} Next we show that if $A$ is a dense subset of a compact
vector space over $\mathbb{F}_2$ and $A+A$ is not much bigger than
$A$ then $2A-2A$ contains a large subspace. \item \emph{(Pullback
and covering)} Finally we use our embedding to pull back this
subspace to a coset in the original setting. A covering argument
then gives us the result.
\end{itemize}
In the remainder of the note we follow through this programme with
our refinement occurring at the second stage.

\subsection{Finding a good model} The appropriate notion of structure preserving was introduced by
Fre{\u\i}man in \cite{GAF}; we record the definition now. If $G$ and
$G'$ are two abelian groups containing the sets $A$ and $A'$
respectively then we say that $\phi:A\rightarrow A'$ is a
\emph{Fre{\u\i}man $s$-homomorphism} if whenever
$a_1,...,a_s,b_1,...,b_s \in A$ satisfy
\begin{equation*}
a_1+...+a_s = b_1+...+b_s
\end{equation*}
we have
\begin{equation*}
\phi(a_1)+...+\phi(a_s) = \phi(b_1)+...+\phi(b_s).
\end{equation*}
If $\phi$ has an inverse which is also an $s$-homomorphism then we
say that $\phi$ is a \emph{Fre{\u\i}man $s$-isomorphism}.

A simple but elegant argument of Green and Ruzsa establishes the
existence of a small vector space into which we can embed our set
via a Fre{\u\i}man isomorphism. Specifically they prove the following
proposition.
\begin{proposition}\label{goodmodel}
\emph{(Proposition 6.1, \cite{BJGIZR})} Suppose that $A$ is a
subset of a vector space over $\mathbb{F}_2$. Suppose that $|A+A|
\leq K|A|$. Then there is a vector space $G'$ over $\mathbb{F}_2$
with $|G'| \leq K^{2s}|A|$, a set $A' \subset G'$, and a Fre{\u\i}man
$s$-isomorphism $\phi:A \rightarrow A'$.
\end{proposition}

\subsection{Bogolyubov's argument} In this section we show that if
$A$ is a subset of a compact vector space over $\mathbb{F}_2$ then
$2A-2A$ contains a large subspace. Originally (in \cite{IZRF})
Ruzsa employed an argument of Bogolubov with the Fourier
transform. This was refined by Chang in \cite{MCC}, and the
improvement of this note rests on a further refinement. We shall
need some notation for the Fourier transform and we record this
now; Rudin, \cite{WR}, includes all the results which we require.

Suppose that $G$ is a compact vector space over $\mathbb{F}_2$.
Write $\wh{G}$ for the dual group, that is the discrete vector
space over $\mathbb{F}_2$ of continuous homomorphisms $\gamma:G
\rightarrow S^1$, where $S^1:=\{z \in \mathbb{C}:|z|=1\}$. $G$ may
be endowed with Haar measure $\mu_G$ normalised so that
$\mu_G(G)=1$ and as a consequence we may define the Fourier
transform $\wh{.}:L^1(G) \rightarrow \ell^\infty(\wh{G})$ which
takes $f \in L^1(G)$ to
\begin{equation*}
\wh{f}: \wh{G} \rightarrow \mathbb{C}; \gamma \mapsto \int_{x \in
G}{f(x)\gamma(-x)d\mu_G(x)}.
\end{equation*}

In \cite{MCC} Chang proved the following result. (Although in
\cite{MCC} it is stated for $\mathbb{Z}/N\mathbb{Z}$, the same
proof applies to any compact abelian group and in particular to
compact vector spaces over $\mathbb{F}_2$.)
\begin{proposition}
Suppose that $G$ is a compact vector space over $\mathbb{F}_2$.
Suppose that $A \subset G$ has density $\alpha$ and $\mu_G(A+A)
\leq K\mu_G(A)$. Then $2A-2A$ contains (up to a null set) a
subspace of codimension $O(K \log \alpha^{-1})$.
\end{proposition}
We prove the following refinement of this.
\begin{proposition}\label{newchang}
Suppose that $G$ is a compact vector space over $\mathbb{F}_2$.
Suppose that $A \subset G$ has density $\alpha$ and $\mu_G(A+A)
\leq K \mu_G(A)$. Then $2A-2A$ contains (up to a null set) a
subspace of codimension $O(K^{1/2}\log \alpha^{-1})$.
\end{proposition}
To prove this we require the following \emph{pure density} version
of the proposition.
\begin{proposition}\label{aaa}
\emph{(Theorem 2.4, \cite{TSASS})} Suppose that $G$ is a compact
vector space over $\mathbb{F}_2$. Suppose that $A \subset G$ has
density $\alpha$. Then $2A-2A$ contains (up to a null set) a
subspace of codimension $O(\alpha^{-1/2})$.
\end{proposition}
The proof in \cite{TSASS} is significantly simpler in the vector
space setting. Since the ideas are important we include the proof
here; the basic technique is iterative.
\begin{lemma}
\emph{(Iteration lemma)}\label{itlem2} Suppose that $G$ is a
compact vector space over $\mathbb{F}_2$. Suppose that $A \subset
G$ has density $\alpha$. Then at least one of the following is
true.
\begin{enumerate}
\renewcommand{\theenumi}{\roman{enumi}}
\item $2A-2A$ contains all of $G$ (up to a null set). \item There
is a subspace $V$ of $\wh{G}$ with dimension 1, an element $x \in
G$ and a set $A' \subset V^\perp$ with the following properties.
\begin{itemize}
\item $x+A' \subset A$; \item $\mu_{V^\perp}(A') \geq
\alpha(1+2^{-1}\alpha^{1/2})$.
\end{itemize}
\end{enumerate}
\end{lemma}
\begin{proof}
As usual with problems of this type studying the sumset $2A-2A$ is
difficult so we turn instead to $g:=\chi_{A} \ast \chi_A \ast
\chi_{-A} \ast \chi_{-A}$ which has support equal to $2A-2A$. One
can easily compute the Fourier transform of $g$ in terms of that
of $\chi_{A}$:
\begin{equation*}
\wh{g}(\gamma)=|\wh{\chi_{A}}(\gamma)|^4 \textrm{ for all } \gamma
\in \wh{G},
\end{equation*}
from which it follows that $g$ is very smooth. Specifically
$\wh{g} \in \ell^{\frac{1}{2}}(\wh{G})$ since
\begin{eqnarray}
\nonumber \sum_{\gamma \in \wh{G}}{|\wh{g}(\gamma)|^{\frac{1}{2}}}
& = & \sum_{\gamma \in
\wh{G}}{|\wh{\chi_{A}}(\gamma)|^2}\\
 \label{smoothness}& = & \alpha
\textrm{ by Parseval's theorem.}
\end{eqnarray}
We may assume that $\mu_G(2A-2A)<1$ since otherwise we are in the
first case of the lemma, so $S:= (2A-2A)^c$ has positive density,
say $\sigma$. Plancherel's theorem gives
\begin{equation*}
0 = \langle \chi_S, g \rangle = \sum_{\gamma \in
\wh{G}}{\overline{\wh{\chi_S}(\gamma)}\wh{g}(\gamma)} \Rightarrow
| \wh{\chi_S}(0_{\wh{G}})\wh{g}(0_{\wh{G}})| \leq \sum_{\gamma
\neq 0_{\wh{G}}}{|\wh{\chi_S}(\gamma)\wh{g}(\gamma)|}.
\end{equation*}
$\wh{g}(0_{\wh{G}})=\alpha^4$, $\wh{\chi_S}(0_{\wh{G}})=\sigma$
and $|\wh{\chi_S}(\gamma)| \leq \|\chi_S\|_1=\sigma$, so the above
yields
\begin{equation*}
\sigma\alpha^4 \leq \sigma\sum_{\gamma \neq
0_{\wh{G}}}{|\wh{g}(\gamma)|} \Rightarrow \alpha^4 \leq
\sum_{\gamma \neq 0_{\wh{G}}}{|\wh{g}(\gamma)|} \textrm{ since
}\sigma>0.
\end{equation*}
Finding a non-trivial character at which $\wh{g}$ is large is now
simple since $\wh{g} \in \ell^{\frac{1}{2}}(\wh{G})$.
\begin{equation*}
\alpha^4 \leq \sup_{\gamma \neq
0_{\wh{G}}}{|\wh{g}(\gamma)|^{\frac{1}{2}}}\left(\sum_{\gamma \in
\wh{G}}{|\wh{g}(\gamma)|^{\frac{1}{2}}}\right) \leq \sup_{\gamma
\neq 0_{\wh{G}}}{|\wh{\chi_A}(\gamma)|^2}. \alpha
\end{equation*}
by (\ref{smoothness}). Rearranging this we have
\begin{equation}\label{usefullinearbias}
\sup_{\gamma \neq 0_{\wh{G}}}{|\wh{\chi_A}(\gamma)|} \geq
\alpha^{\frac{3}{2}}.
\end{equation}
The set $\Gamma:=\{\gamma \in \wh{G}: |\wh{\chi_A}(\gamma)| \geq
\alpha^{\frac{3}{2}} \}$ has size at most $\alpha^{-2}$ since
\begin{equation*}
|\Gamma|.(\alpha^{\frac{3}{2}})^{2} \leq \sum_{\gamma \in
\wh{G}}{|\wh{\chi_A}(\gamma)|^{2}} = \alpha \textrm{ by Parseval's
theorem.}
\end{equation*}
It follows that the supremum in (\ref{usefullinearbias}) is really
a maximum and we may pick a character $\gamma$ which attains this
maximum. We now proceed with a standard
$L^\infty$-density-increment argument. Let
$V:=\{0_{\wh{G}},\gamma\}$ and $f:=\chi_{A} - \alpha$. Then
\begin{equation*}
\int{f \ast \mu_{V^\perp} d\mu_G}=0 \textrm{ and } \|f \ast
\mu_{V^\perp}\|_1 \geq \|\wh{f}\wh{\mu_{V^\perp}}\|_\infty =
|\wh{\chi_{A}}(\gamma)|.
\end{equation*}
Adding these we conclude that
\begin{eqnarray*}
|\wh{\chi_{A}}(\gamma)| & \leq &2\int{(f \ast
\mu_{V^\perp})_+d\mu_G}\\ & = & 2\int{(\chi_{A} \ast \mu_{V^\perp}
- \alpha)_+d\mu_G}\\ & \leq & 2(\|\chi_{A} \ast
\mu_{V^\perp}\|_\infty - \alpha).
\end{eqnarray*}
$\chi_A \ast \mu_{V^\perp}$ is continuous so there is some $x \in
G$ with
\begin{equation*}
\chi_A \ast \mu_{V^\perp}(x) = \|\chi_A \ast
\mu_{V^\perp}\|_\infty \geq \alpha(1+2^{-1}\alpha^{1/2}).
\end{equation*}
The result follows on taking $A'=x+A$.
\end{proof}
\begin{proof}[Proof of Proposition \ref{aaa}.]
We define a nested sequence of finite dimensional subspaces $V_0
\leq V_1 \leq ... \leq \wh{G}$, elements $x_k \in V_k^\perp$ and
subsets $A_k$ of $V_k^\perp$ with density $\alpha_k$, such that
$x_k+ A_{k} \subset A_{k-1}$. We begin the iteration with
$V_0:=\{0_{\wh{G}}\}$, $A_0:=A$ and $x_0=0_{G}$.

Suppose that we are at stage $k$ of the iteration. If
$\mu_{V_k^\perp}(2A_k-2A_k)<1$ then we apply Lemma \ref{itlem2} to
$A_k$ considered as a subset of $V_k^\perp$. We get a vector space
$V_{k+1}$ with $\dim V_{k+1} = 1+ \dim V_k$, an element $x_{k+1}
\in G$ and a set $A_{k+1}$ such that
\begin{equation*}
x_{k+1} + A_{k+1} \subset A_k \textrm{ and } \alpha_{k+1} \geq
\alpha_k(1+2^{-1}\alpha_k^{1/2}).
\end{equation*}
It follows from the density increment that if
$m_k=2\alpha_k^{-1/2}$ then $\alpha_{k+m_k} \geq 2\alpha_k$.
Define the sequence $(N_{l})_l$ recursively by $N_0=0$ and
$N_{l+1}=m_{N_l}+N_l$. The density $\alpha_{N_l}$ is easily
estimated:
\begin{equation*}
\alpha_{N_l} \geq 2^l\alpha \textrm{ and } N_l \leq
\sum_{s=0}^l{2\alpha_{N_s}^{-1/2}} \leq 2
\alpha^{-1/2}\sum_{s=0}^l{2^{-s/2}} = O(\alpha^{-1/2}).
\end{equation*}
Since density cannot be greater than 1 there is some stage $k$
with $k = O(\alpha^{-1/2})$ when the iteration cannot proceed i.e.
for which $2A_k-2A_k$ contains all of $V_k^\perp$ (except for a
null set). By construction of the $A_k$s there is a translate of
$A_k$ which is contained in $A_0=A$ and hence $2A_k-2A_k$ is
contained in $2A-2A$. It follows that $2A-2A$ contains (up to a
null set) a subspace of $G$ of codimension $k=O(\alpha^{-1/2})$.
\end{proof}
The key ingredient in the proof of Proposition \ref{newchang} is
the following iteration lemma, which has a number of similarities
with Lemma \ref{itlem2}.
\begin{lemma}\label{itl}
Suppose that $G$ is a compact vector space over $\mathbb{F}_2$.
Suppose that $A, B \subset G$ have $\mu_G(A+B) \leq K \mu_G(B)$.
Write $\alpha$ for the density of $A$. Then at least one of the
following is true.
\begin{enumerate}
\renewcommand{\theenumi}{\roman{enumi}}
\item $B$ contains (up to a null set) a subspace of codimension
$O(K^{1/2})$. \item There is a subspace $V$ of $\wh{G}$ with
dimension 1, elements $x,y \in G$ and sets $A',B' \subset V^\perp$
with the following properties.
\begin{itemize}
\item $x+A' \subset A \textrm{ and } y+B' \subset B$; \item
$\mu_{V^\perp}(A') \geq \alpha(1+2^{-3/2}K^{-1/2})$; \item
$\mu_{V^\perp}(A'+B') \leq K\mu_{V^\perp}(B')$.
\end{itemize}
\end{enumerate}
\end{lemma}
\begin{proof}
If $\mu_G(B) \geq (2K)^{-1}$ then we apply Proposition \ref{aaa}
to get that $B$ contains (up to a null set) a subspace of
codimension $O(K^{1/2})$ and we are in the first case of the
lemma. Hence we assume that $\mu_G(B) \leq (2K)^{-1}$.

Write $\beta$ for the density of $B$. We have
\begin{eqnarray}
\nonumber (\alpha\beta)^2 & = & \left(\int{\chi_A \ast
\chi_Bd\mu_G}\right)^2\\ \nonumber & \leq & \mu_G(A+B)\int{(\chi_A
\ast \chi_B)^2d\mu_G} \textrm{ by Cauchy-Schwarz,}\\ \nonumber  &
\leq & K \beta \int{(\chi_A \ast \chi_B)^2d\mu_G} \textrm{ by
hypothesis,}\\ \label{ytu} & = & K\beta \sum_{\gamma \in
\wh{G}}{|\wh{\chi_A}(\gamma)|^2|\wh{\chi_B}(\gamma)|^2} \textrm{
by Parseval's theorem.}
\end{eqnarray}
The main term in the sum on the right is the contribution from the
trivial character, in particular
\begin{equation*}
|\wh{\chi_A}(0_{\wh{G}})|^2|\wh{\chi_B}(0_{\wh{G}})|^2=\alpha^2\beta^2,
\end{equation*}
while
\begin{eqnarray*}
\sum_{\gamma \neq
0_{\wh{G}}}{|\wh{\chi_A}(\gamma)|^2|\wh{\chi_B}(\gamma)|^2} & \leq
& \sup_{\gamma \neq 0_{\wh{G}}}{|\wh{\chi_A}(\gamma)|^2}
\sum_{\gamma \in \wh{G}}{|\wh{\chi_B}(\gamma)|^2}\\ & = & \beta
\sup_{\gamma \neq 0_{\wh{G}}}{|\wh{\chi_A}(\gamma)|^2}\textrm{ by
Parseval's theorem for $\chi_B$}.
\end{eqnarray*}
Putting these last two observations in (\ref{ytu}) gives
\begin{equation*}
\alpha^2\beta^2 \leq K\beta^3\alpha^2 + K\beta^2\sup_{\gamma \neq
0_{\wh{G}}}{|\wh{\chi_A}(\gamma)|^2}.
\end{equation*}
Since $K\beta \leq 2^{-1}$ we can rearrange this to conclude that
\begin{equation}\label{usefullinearbias2}
\sup_{\gamma \neq 0_{\wh{G}}}{|\wh{\chi_A}(\gamma)|} \geq
(2K)^{-1/2}\alpha.
\end{equation}
The set $\Gamma:=\{\gamma \in \wh{G}: |\wh{\chi_A}(\gamma)| \geq
(2K)^{-1/2}\alpha \}$ has size at most $2K\alpha^{-1}$ since
\begin{equation*}
|\Gamma|.((2K)^{-1/2}\alpha)^{2} \leq \sum_{\gamma \in
\wh{G}}{|\wh{\chi_A}(\gamma)|^{2}} = \alpha \textrm{ by Parseval's
theorem.}
\end{equation*}
It follows that the supremum in (\ref{usefullinearbias2}) is
really a maximum and we may pick a character $\gamma$ which
attains this maximum. We now proceed with a standard
$L^\infty$-density-increment argument. Let
$V:=\{0_{\wh{G}},\gamma\}$ and $f:=\chi_{A} - \alpha$. Then
\begin{equation*}
\int{f \ast \mu_{V^\perp} d\mu_G}=0 \textrm{ and } \|f \ast
\mu_{V^\perp}\|_1 \geq \|\wh{f}\wh{\mu_{V^\perp}}\|_\infty =
|\wh{\chi_{A}}(\gamma)|.
\end{equation*}
Adding these we conclude that
\begin{eqnarray*}
|\wh{\chi_{A}}(\gamma)| & \leq &2\int{(f \ast
\mu_{V^\perp})_+d\mu_G}\\ & = & 2\int{(\chi_{A} \ast \mu_{V^\perp}
- \alpha)_+d\mu_G}\\ & \leq & 2(\|\chi_{A} \ast
\mu_{V^\perp}\|_\infty - \alpha).
\end{eqnarray*}
Since $\chi_A \ast \mu_{V^\perp}$ is continuous it follows that
there is some $x$ for which
\begin{equation*}
\chi_A \ast \mu_{V^\perp}(x) \geq \alpha (1+2^{-3/2}K^{-1/2}).
\end{equation*}
Let $x'+V^\perp:=G \setminus (x+V^\perp)$ be the \emph{other}
coset of $V^\perp$ in $G$. Write $A_1=A \cap (x+V^\perp)$,$B_1=B
\cap (x+V^\perp)$ and $B_2=B \cap (x'+V^\perp)$. Now $A_1 \subset
A$ so
\begin{equation*}
(A_1+B_1) \cup (A_1+B_2) \subset A + B_1 \cup B_2,
\end{equation*}
and $A_1+B_1 \subset V^\perp$ while $A_1+B_2 \subset x+x'+V^\perp$
so these two sets are disjoint and we conclude that
\begin{eqnarray*}
\mu_G(A_1+B_1) + \mu_G(A_1+B_2) & = & \mu_G( (A_1+B_1) \cup
(A_1+B_2))
\\ & \leq & \mu_G( A + B_1 \cup B_2)\\ & \leq & K
\mu_G(B_1 \cup B_2)\textrm{ by hypothesis }\\ & \leq &
K(\mu_G(B_1) + \mu_G(B_2)).
\end{eqnarray*}
Hence, by averaging, there is some $i$ such that
\begin{equation*}
\mu_G(A_1+B_i) \leq K\mu_G(B_i).
\end{equation*}
We take $A'=x+A_1$ and, if $i=1$, $B'=x+B_1$ and $y=x$, while if
$i=2$, $B'=x'+B_2$ and $y=x'$. The result follows.
\end{proof}

\begin{proof}[Proof of Proposition \ref{newchang}.]
We define a nested sequence of finite dimensional subspaces $V_0
\leq V_1 \leq ... \leq \wh{G}$, elements $x_k,y_k \in V_k^\perp$,
and subsets $A_k$ and $B_k$ of $V_k^\perp$ such that $A_{k} + x_k
\subset A_{k-1}$ and $B_k + y_k \subset B_{k-1}$ and
$\mu_{V_k^\perp}(A_k+B_k) \leq K \mu_{V_k^\perp}(B_k)$. We write
$\alpha_k$ for the density of $A_k$ in $V_k^\perp$. Begin the
iteration with $V_0:=\{0_{\wh{G}}\}$, $B_0=A_0:=A$ and
$x_0=y_0=0_{G}$.

Suppose that we are at stage $k$ of the iteration. We apply Lemma
\ref{itl} to $A_k$ and $B_k$ inside $V_k^\perp$ (which we can do
since $\mu_{V_k^\perp}(A_k+B_k) \leq K \mu_{V_k^\perp}(B_k)$). It
follows that either $2B_k-2B_k$ contains (up to a null set) a
subspace of codimension $O(K^{1/2})$ in $V_k^\perp$ or we get a
subspace $V_{k+1} \leq \wh{V_{k}}$ with $\dim V_{k+1}=1+\dim V_k$,
elements $x_{k+1},y_{k+1} \in V_k^\perp$ and sets $A_{k+1}$ and
$B_{k+1}$ with the following properties.
\begin{itemize}
\item $x_{k+1}+A_{k+1} \subset A_k \textrm{ and } y_{k+1}+B_{k+1}
\subset B_k$; \item $\mu_{V^\perp}(A_{k+1}) \geq
\alpha_{k}(1+2^{-3/2}K^{-1/2})$; \item
$\mu_{V^\perp}(A_{k+1}+B_{k+1}) \leq K\mu_{V^\perp}(B_{k+1})$.
\end{itemize}
It follows from the density increment that if $m=2^{3/2}K^{1/2}$
then $\alpha_{k+m} \geq 2\alpha_k$, and hence the iteration must
terminate (because density can be at most 1) at some stage $k$
with $k=O(K^{1/2}\log \alpha^{-1})$. The iteration terminates if
$2B_k-2B_k$ contains (up to a null set) a subspace of codimension
$O(K^{1/2})$ in $V_k^\perp$, from which it follows that $2A-2A
\supset 2B_k -2B_k$ contains (up to a null set) a subspace of
codimension $k+O(K^{1/2}) = O(K^{1/2}\log \alpha^{-1})$.
\end{proof}

\subsection{Pullback and covering} We now complete the proof of
the main theorem using a covering argument.

We are given $A \subset G$ finite with $|A+A| \leq K|A|$. By
Proposition \ref{goodmodel} there is a finite vector space $G'$
with $|G'| \leq K^{16}|A|$ and a subset $A'$ with $A'$ Fre{\u\i}man
8-isomorphic to $A$. It follows that
\begin{equation*}
\mu_{G'}(A') \geq K^{-16} \textrm{ and } \mu_{G'}(A'+A') \leq
K\mu_{G'}(A').
\end{equation*}
We apply Proposition \ref{newchang} to conclude that $2A'-2A'$
contains a subspace of codimension $O(K^{1/2}\log K)$. However,
$A$ is 8-isomorphic to $A'$ so $2A-2A$ is 2-isomorphic to
$2A'-2A'$ and it is easy to check that the 2-isomorphic pullback
of a subspace is a coset so $2A-2A$ contains a coset of size
\begin{equation*}
2^{-O(K^{1/2}\log K)}|G'| \geq 2^{-O(K^{1/2}\log K)}|A|.
\end{equation*}

The following covering result of Chang \cite{MCC} converts this
large coset contained in $2A-2A$ into a small coset containing
$A$. It is true in more generality than we state; we only require
the version below.
\begin{proposition}
Suppose that $G$ is a vector space over $\mathbb{F}_2$. Suppose
that $A \subset G$ is a finite set with $|A+A| \leq K|A|$. Suppose
that $2A-2A$ contains a coset of size $\eta|A|$. Then $A$ is
contained in a coset of size at most $2^{O(K\log K\eta^{-1})}|A|$.
\end{proposition}
Theorem \ref{maintheorem} follows immediately from this
proposition and the argument preceding it.

\section{Concluding remarks}

It is worth making a couple of concluding remarks. First, all the
implicit constants in the work are effective however they are not
particularly neat or significant so it does not seem to be
important to calculate them. Secondly, and more importantly, it
seems likely that one could modify Proposition \ref{newchang} to
fall within the more general framework of approximate groups
pioneered by Bourgain in \cite{JB}. It does not seem that this
would lead to any improvement in Fre{\u\i}man's theorem for
$\mathbb{Z}$, essentially because of the need to narrow the Bohr
sets at each stage of the iteration.

\section*{Acknowledgements} I should like to thank Tim Gowers and
Ben Green for encouragement and supervision.

\bibliographystyle{alpha}

\bibliography{master}

\end{document}